\documentclass[a4paper,12pt]{amsart}
\usepackage{amsmath,latexsym,amscd}

\newtheorem{thm}{Theorem}[section]
\newtheorem{cor}[thm]{Corollary}
\newtheorem{lemma}[thm]{Lemma}
\newtheorem{prop}[thm]{Proposition}
\theoremstyle{definition}
\newtheorem{defn}[thm]{Definition}
\theoremstyle{remark}
\newtheorem{remark}[thm]{Remark}

\newtheorem{setup}[thm]{}

\newtheorem{claim}[thm]{Claim}

\newcommand{\bQ}{{\Bbb Q}}
\newcommand{\roundup}[1]{\ulcorner{#1}\urcorner}

\newcommand{\fei}[1]{\phi_{#1}}
\newcommand{\Fi}[1]{\Phi_{#1}}
\newcommand{\lrw}{\longrightarrow}
\newcommand{\simlin}{\sim_{\text{lin}}}
\newcommand{\simnum}{\sim_{\text{num}}}
\newcommand{\Co}[1]{{\mathcal O}_{#1}}

\begin{document}

\title
{On the $\bQ$-divisor method and its application}
\author{M. Chen}
\date{}
\address{\rm
Institute of Mathematics, Fudan University, Shanghai, 200433, PR
China} \email{mchen@fudan.edu.cn}
\address{and}
\address{\rm The Institute of Mathematical Sciences, The Chinese
University of Hong Kong}

\thanks{This paper is supported by The Institute of Mathematical
Sciences, The Chinese University of Hong Kong. The project is also
supported by the National Natural Science Foundation of China
(No.10131010), Shanghai Scientific $\&$ Technical Commission
(Grant 01QA14042) and SRF for ROCS, SEM}

\maketitle
\begin{abstract}
For a smooth projective 3-fold of general type, we prove that the
relative canonical stability $\mu_s(3)\le 8$. This is induced by
our improved result of Koll\'ar: the m-canonical map of $X$ is
birational onto its image whenever $m\ge 5k+6$, provided
$P_k(X)\ge 2$. The ${\mathbb Q}$-divisor method is intensively
developed to prove our results.
\end{abstract}

\pagestyle{myheadings} \markboth{\hfill M. Chen\hfill}{\hfill The
Q-divisor method\hfill}

\section*{\rm Introduction}

Let $X$ be a smooth projective threefold with big $K_X$. Many
authors have studied the pluricanonical systems or the
pluricanonical maps of $X$(see \cite{Ben}, \cite{Ch1}, \cite{E-L},
\cite{Ka1}, \cite{Kol}, \cite{L1}, \cite{L2}, \cite{Ma},
\cite{Sho}, \cite{Wi} etc). Suppose $h^0(X, kK_X)\ge 2$. J.
Koll\'ar(\cite{Kol}, Corollary 4.8) first proved that the
$(11k+5)$-canonical map is birational. Then, in \cite{Ch2}, it has
been proved that either the $(7k+3)$-canonical map or the
$(7k+5)$-canonical map is birational. We denote by $\fei{m}$ the
m-canonical map of $X$. Since it's still unclear whether the
birationality of  $\fei{m}$ does imply the birationality of
$\fei{m+1}$, we hope to find certain "stable property" of
$\fei{m}$. In this paper, we modify the ${\Bbb Q}$-divisor method
and then present much better results than in \cite{Kol} and
\cite{Ch2,Ch3}. The ${\Bbb Q}$-divisor method was originally
developed by Kawamata, Reid, Shokorov and others in connection
with the minimal model program initiated by Mori. It was also
exploited much effectively, for various purpose, by many authors
such as Ein, Koll\'ar, Lazarsfeld, Viehweg and so on. As far as
our method can tell here, the results are as follows.

\begin{thm}\label{T:0.1}  Let $X$ be a smooth projective 3-fold of general type.
Suppose the $k$-th plurigenus $P_k(X)\ge 2$. Then the $m$-canonical map is
birational onto its image for all $m\ge 5k+6$.
\end{thm}

\begin{thm}\label{T:0.2}  Let $X$ be a smooth projective 3-fold of general type.
Suppose $p_g(X)\ge 2$. Then $\fei{8}$ is birational onto its image.
\end{thm}

The base field is always supposed to be algebraically closed of characteristic 0.

For readers' convenience, let us recall the definition of the so-called {\it relative canonical stability} as follows.

Let $X$ be a nonsingular projective variety of general type of dimension
$n$. We define

$k_0(X):=min\{\ k|\ P_k(X)\ge 2\};$

$k_s(X):=min\{\ k|\ \phi_m$ is birational for $m\ge k\}$, which is called {\it
the canonical stability of} $X$;

$\mu_s(X):=\frac{k_s(X)}{k_0(X)}$, which is called {\it the relative
canonical stability} of $X$.

$\mu_s(n):=\text{sup}\{\mu_s(X)| X$ is a smooth projective n-fold of general type$\}$,
which is referred to as {\it the n-th relative canonical stability}.

It's well-known that $\mu_s(1)=3$ and $\mu_s(2)=5$. A conjecture
concerning $\mu_s$ is that $\mu_s(n)=2n+1$. Theorem \ref{T:0.1}
and Theorem \ref{T:0.2} imply the following

\begin{cor} $\mu_s(3)\le 8$.
\end{cor}

\begin{remark} $\mu_s(n)$ is very interesting for $n\ge 3$. For example, if $\mu_s(4)<+\infty$,
then the following is true:

there is a constant $m_0$ such that $\fei{m_0}$ is birational for all smooth projective 3-fold $X$ of general type with
$p_g(X)>0$.

\noindent In fact, taking the product of $X$ with a curve of genus $\ge 2$, one can realize the above statement very easily.
Unfortunately, few information found on $\mu_s(4)$ yet.
\end{remark}

\section{\rm Key lemmas}\label{S:1}

Let $S$ be a smooth projective surface and $D$ a divisor on $S$. We would
like to know when the system $|K_S+D|$ gives a birational rational map onto
its image. It is well-known that Reider (\cite{Rdr}) gave much effective
results when $D$ is nef and big. Because, in our cases, $D$ may be not nef, we
hope to find an effective sufficient condition in order to treat our situation.

Suppose $M$ is a divisor on $S$ with $h^0(S, M)\ge 2$ and $|M|$ is
movable. Taking necessary blow-ups $\pi:S'\lrw S$, along the
indeterminency of the system, such that the movable part of
$|\pi^*(M)|$ is basepoint free, we can get a morphism
$g:=\Fi{M}\circ \pi: S'\lrw W$ where $W$ is the image of $S'$ in
${\Bbb P}^{h^0(S, M)-1}$. If $\dim(W)=1$, we can take the Stein
factorization $g: S'\overset{f}\lrw B\lrw W$, where $B$ is a
smooth projective curve and $f$ is a fibration. Denote $b:=g(B)$.

\begin{defn}\label{D:1.1} If $b=0$, we say that $|M|$ is {\it composed of a
rational pencil of curves}. Otherwise, $|M|$ is {\it composed of an
irrational pencil of curves.}
\end{defn}

We can write $\pi^*(M)\simlin M_0+Z_0$ where $M_0$ is the movable part and
$Z_0$ the fixed one. According to Bertini's theorem, a general member
$C\in |M_0|$ is a smooth irreducible curve if $|M_0|$ is not composed of
pencil of curves. Otherwise, we can write $M_0\simlin \sum_{i=1}^mC_i$ where
the $C_i$'s are fibers of the fibration $f$ for all $i$. Denote by $C$ a
general fiber of $f$ in the latter case. Then $C$ is still an irreducible
smooth curve.

\begin{defn}\label{D:1.2} In the above setting, $C$ is called
{\it a generic irreducible element of} $|M_0|$.  Meanwhile, $\pi(C)$ is
called {\it a generic irreducible element of} $|M|$. Note that, however,
$\pi(C)$ may be a singular curve.
\end{defn}

Now we can set up the following lemma.

\begin{lemma}\label{L:1.3} Let $S$ be a smooth projective surface and $D$ a divisor on $S$. The rational map $\Fi{K_S+D}$ is birational onto its image if the following conditions hold.

(i) There is a divisor $M$ on $S$ with $h^0(S, M)\ge 2$ and $|M|$ movable such that
$K_S+D\ge M$. If $|M|$ is composed of an irrational pencil of curves, denoting by $C$ its generic irreducible element, $D-2C\ge \roundup{A_0}$ holds
for certain nef and big ${\Bbb Q}$-divisor $A_0$ on $S$ (with $A_0\cdot C>1$ in the case
$\kappa(S)=-\infty$).

(ii) $D-C\ge\roundup{A}$ holds for certain nef and big ${\Bbb Q}$-divisor $A$ on $S$ with $A\cdot C>2$, where $C$ is a generic irreducible element of $|M|$.
\end{lemma}
\begin{proof}
It's clear that one may suppose that $|M|$ is base point free. Taking the
Stein-factorization of $\Fi{M}$, we can get
$$\Fi{M}: X\overset{f}{\lrw} W\overset{s}{\lrw} {\Bbb P}^{h^0(M)-1},$$
where $f$ has connected fibers. A generic irreducible element of $|M|$ is a smooth
projective curve.

If $|M|$ is not composed of a pencil of curves, then it's sufficient to verify
the birationality of $\Fi{K_S+D}|_{C}$ by virtue of Tankeev's principle (\cite{Ta}, Lemma 2).

If $|M|$ is composed of a pencil of curves, then $f$ is a fibration onto the
smooth curve $W$. In this case, $C$ is a general fiber of $f$.
When $W$ is a rational curve, $|K_S+D|$ can distinguish different
fibers of $f$ since $K_S+D\ge C$. When $W$ is irrational, suppose $C_1$ and $C_2$
are two general fibers of $f$. By assumption, we have
$$D-C_1-C_2=\roundup{A_0}+F,$$
where $F$ is an effective divisor on $S$ and $A_0$ is a nef and big
${\Bbb Q}$-divisor. We consider the system $|K_S+D-F|$. According to
Kawamata-Viehweg vanishing theorem, we have the surjective map
$$H^0(S, K_S+D-F)\lrw H^0(C_1, K_{C_1}+G_1)\oplus H^0(C_2, K_{C_2}+G_2).$$
By (i), we have $h^0(C_1, K_{C_1}+G_1)>0$ and $h^0(C_2, K_{C_2}+G_2)>0$. Thus
$|K_S+D-F|$ can distinguish $C_1$ and $C_2$, so can $|K_S+D|$. Therefore
it's also sufficient to verify the birationality of $\Fi{K_S+D}|_{C}$.

By (ii), we can write $D-C=\roundup{A}+F_1$ where $F_1$ is an effective divisor and
$A$ is a nef and big ${\Bbb Q}$-divisor. We consider the system $|K_S+D-F_1|$.
By vanishing theorem, we have the surjective map
$$H^0(S, K_S+D-F_1)\lrw H^0(C, K_C+G)$$
where $G$ is a divisor on $C$ with $\deg (G)\ge 3$. Thus
$\Fi{K_S+D-F_1}|_{C}$ is an embedding, so is $\Fi{K_S+D}|_{C}$.
\end{proof}

\begin{lemma}\label{L:1.4}
Let $S$ be a smooth projective surface of general type. Then
$K_S+\roundup{A}+D$ is effective if $A$ is a nef and big ${\Bbb Q}$-divisor
and if $h^0(S, D)\ge 2$.
\end{lemma}
\begin{proof}
We may suppose that
$|D|$ is basepoint free. Denote by $C$ a generic irreducible element of
$|D|$. Then the vanishing theorem gives the exact sequence
$$H^0(S, K_S+\roundup{A}+C)\lrw H^0(C, K_C+H)\lrw 0,$$
where $H:=\roundup{A}|_C$ is a divisor of positive degree. It is obvious that
$h^0(C, K_C+H)\ge 2$ since $C$ is a curve of genus $\ge 2$. The proof is completed.
\end{proof}

\begin{lemma}\label{L:1.5} (\cite{Ch3}, Lemma 2.7)
Let $X$ be a smooth projective variety of dimension $\ge 2$. Let $D$ be a divisor
on $X$, $h^0(X, \Co{X}(D))\ge 2$ and $S$ be a smooth irreducible divisor on $X$ such that $S$ is not a fixed component of $|D|$. Denote by $M$ the movable part of $|D|$ and by $N$ the movable part of $|D|_S|$ on $S$. Suppose the natural restriction map
$$H^0(X, \Co{X}(D))\overset\theta\lrw H^0(S, \Co{S}(D|_S))$$
is surjective. Then $M|_S\ge N$ and thus
$$h^0(S,\Co{S}(M|_S))=h^0(S,\Co{S}(N))=h^0(S,\Co{S}(D|_S)).$$
\end{lemma}

\section{\rm The case with $k_0\ge 2$}\label{S:2}

Sometimes, for simplicity, we denote $k_0(X)$ and $k_s(X)$ by $k_0$ and
$k_s$  respectively.

\begin{prop}\label{P:2.1} Let $X$ be a minimal projective 3-fold of general
type with only ${\Bbb Q}$-factorial terminal singularities. If
$\dim\fei{k_0}(X)\ge 2$, then $P_m(X)\ge 2$ for all $m\ge 2k_0$.
\end{prop}
\begin{proof}
First we take a birational modification $\pi: X'\lrw X$, according to
Hironaka, such that

(1) $X'$ is smooth;

(2) the movable part of $|k_0K_{X'}|$ defines a morphism;

(3) the fractional part of $\pi^*(K_X)$ has supports with only normal
crossings.

Denote by $S_0:=S_{k_0}$ a generic irreducible element of the movable part
of $|k_0K_{X'}|$. Then $S_0$ is a smooth projective surface of general type
by Bertini's theorem. By the vanishing theorem, we have the exact sequence
\begin{align*}
&H^0(X', K_{X'}+\roundup{(t+k_0)\pi^*(K_X)}+S_0)\\
\lrw &H^0(S_0, K_{S_0}+
\roundup{(t+k_0)\pi^*(K_X)}|_{S_0})\lrw 0,
\end{align*}
where $t\ge 0$ is a given integer and
$$\roundup{(t+k_0)\pi^*(K_X)}|_{S_0}\ge \roundup{t\pi^*(K_X)|_{S_0}}+D_0,$$
$D_0:=S_0|_{S_0}$ has the property
$h^0(S_0,D_0)\ge 2$ according to the assumption.

If $t=0$, then
$$P_{2k_0+1}(X)\ge h^0(S_0, K_{S_0}+D_0)\ge 2$$
 by Lemma 1.2 of \cite{Ch2}.

If $t>0$, we still have the following exact sequence
$$H^0(S_0, K_{S_0}+\roundup{t\pi^*(K_X)|_{S_0}}+C)\lrw
H^0(K_C+G)\lrw 0,$$
where $C$ is a generic irreducible element of the movable part of $|D_0|$ and
$G:=\roundup{t\pi^*(K_X)|_{S_0}}|_C$
is a divisor of positive degree on $C$.
Since $C$ is a curve of genus $\ge 2$, we have
$h^0(C,K_C+G)\ge 2.$
Thus we can  see that
$P_{2k_0+t+1}\ge 2$. The proof is completed.
\end{proof}

\begin{cor}\label{C:2.2} Let $X$ be a minimal projective 3-fold of general
type with only ${\Bbb Q}$-factorial terminal singularities. If
$\fei{k_0}$ is birational, then $k_s\le 3k_0$.
\end{cor}
\begin{proof} This is obvious according to Proposition 2.1.
\end{proof}

\begin{thm}\label{T:2.3}  Let $X$ be a minimal projective 3-fold of general
type with only ${\Bbb Q}$-factorial terminal singularities.
If $\dim\fei{k_0}(X)=3$ and $\fei{k_0}$ is not birational, then
$k_s\le 3k_0+2$.
\end{thm}
\begin{proof}
Taking the same modification $\pi:X'\lrw X$ as in the proof of Proposition 2.1,
we still denote by $S_0$ the general member of the movable part of
$|k_0K_{X'}|$. Note that both $|k_0K_{X'}|$ and $|\roundup{k_0\pi^*(K_X)}|$
have the same movable part. For a given integer $t>0$, we have
$$K_{X'}+\roundup{(t+2k_0)\pi^*(K_X)}+S_0\le (t+3k_0+1)K_{X'}.$$
It is sufficient to prove the birationality of the rational map defined by
$$|K_{X'}+\roundup{(t+2k_0)\pi^*(K_X)}+S_0|.$$
Because
$K_{X'}+\roundup{(t+2k_0)\pi^*(K_X)}$
is effective by the proof of Proposition 2.1, we have
$$K_{X'}+\roundup{(t+2k_0)\pi^*(K_X)}+S_0\ge S_0.$$
Thus we only need to prove the birationality of
$$\Fi{K_{X'}+\roundup{(t+2k_0)\pi^*(K_X)}+S_0}\bigm|_{S_0}.$$
We have the following exact sequence by the vanishing theorem
\begin{align*}
&H^0(X, K_{X'}+\roundup{(t+2k_0)\pi^*(K_X)}+S_0)\\
\lrw
&H^0(S_0, K_{S_0}+\roundup{(t+2k_0)\pi^*(K_X)}|_{S_0})\lrw 0,
\end{align*}
which means
$$|K_{X'}+\roundup{(t+2k_0)\pi^*(K_X)}+S_0|\bigm|_{S_0}=
|K_{S_0}+\roundup{(t+2k_0)\pi^*(K_X)}|_{S_0}|.$$
Noting that
$$K_{S_0}+\roundup{(t+2k_0)\pi^*(K_X)}|_{S_0}\ge
K_{S_0}+\roundup{t\pi^*(K_X)|_{S_0}}+2L_0,$$
where $L_0:=S_0|_{S_0}$, we want to show that
$\Fi{K_{S_0}+\roundup{t\pi^(K_X)|_{S_0}}+2L_0}$
is birational. Because $|L_0|$ gives a generically finite map, we see from
Lemma 1.4 that
$K_{S_0}+\roundup{t\pi^*(K_X)|_{S_0}}+L_0$
is effective. On the other hand,
let $C$ be a generic irreducible element of $|L_0|$, then $\dim\Fi{L_0}(C)=1$.
So $L_0\cdot C\ge 2$ and thus $(t\pi^*(K_X)|_{S_0}+L_0)\cdot C>2$. By Lemma 1.3,
$|K_{S_0}+\roundup{t\pi^*(K_X)|_{S_0}}+2L_0|$
gives a birational map. The proof is completed.
\end{proof}

\begin{thm}\label{T:2.4}  Let $X$ be a minimal projective 3-fold of general
type with only ${\Bbb Q}$-factorial terminal singularities.
If $\dim\fei{k_0}(X)=2$, then $k_s\le 4k_0+4$.
\end{thm}
\begin{proof}
First we take the same modification $\pi:X'\lrw X$ as in the proof of
Proposition 2.1. We also suppose that $S_0$ is the movable part of $|k_0K_{X'}|$.
For a given integer $t>0$, we obviously have
$$K_{X'}+\roundup{(t+2k_0+2)\pi^*(K_X)}+2S_0\le (t+4k_0+3)K_{X'}.$$
Thus it is sufficient to verify the birationality of the rational map defined
by
$$|K_{X'}+\roundup{(t+2k_0+2)\pi^*(K_X)}+2S_0|.$$
By Proposition 2.1,
$$K_{X'}+\roundup{(t+2k_0+2)\pi^*(K_X)}+S_0$$
is effective.
Thus we only have to prove the birationality of the restriction
$$\Fi{K_{X'}+\roundup{(t+2k_0+2)\pi^*(K_X)}+2S_0}\bigm|_{S_0}$$
for the general $S_0$. The vanishing theorem gives the exact sequence
\begin{align*}
&H^0(X',K_{X'}+\roundup{(t+2k_0+2)\pi^*(K_X)}+2S_0)\\
\lrw
&H^0(S_0, K_{S_0}+\roundup{(t+2k_0+2)\pi^*(K_X)}\bigm|_{S_0}+S_0|_{S_0})\lrw
0.
\end{align*}
This means
$$\Fi{K_{X'}+\roundup{(t+2k_0+2)\pi^*(K_X)}+2S_0}\bigm|_{S_0}
=\Fi{K_{S_0}+\roundup{(t+2k_0+2)\pi^*(K_X)}|_{S_0}+S_0|_{S_0}}.$$
Suppose $M_{2k_0+2}$ is the movable part of $|(2k_0+2)K_{X'}|$. We
have to study some property of $\bigm|M_{2k_0+2}|_{S_0}\bigm|$.
Note that $M_{2k_0+2}$ is also the movable part of
$$|\roundup{(2k_0+2)\pi^*(K_X)}|.$$
 We have
$K_{X'}+\roundup{\pi^*(K_X)}+2S_0\le (2k_0+2)K_{X'}.$
The vanishing theorem gives the exact sequence
$$H^0(X',K_{X'}+\roundup{\pi^*(K_X)}+2S_0)\overset{\alpha}\lrw
H^0(S_0, K_{S_0}+\roundup{\pi^*(K_X)}|_{S_0}+L_0)\lrw 0,$$
where $L_0:=S_0|_{S_0}$. Denote by $M_{2k_0+2}'$ the movable part of
$|K_{X'}+\roundup{\pi^*(K_X)}+2S_0|$
and by $G$ the movable part of
$|K_{S_0}+\roundup{\pi^*(K_X)}|_{S_0}+L_0|.$
By Lemma 2.3, we have
$G\le M_{2k_0+2}'|_{S_0}\le M_{2k_0+2}|_{S_0}.$
Noting that $|L_0|$ is a free pencil, we can suppose $C$ is a generic irreducible
element of $|L_0|$. Now the key step is to show that $\dim\Fi{G}(C)=1$.
In fact, the vanishing theorem gives
$$|K_{S_0}+\roundup{\pi^*(K_X)|_{S_0}}+L_0|\bigm|_C=|K_C+D|,$$
where $D:=\roundup{\pi^*(K_X)|_{S_0}}\bigm|_C$ is a divisor of positive degree.
Because $C$ is a curve of genus $\ge 2$, $|K_C+D|$ gives a finite map.
This shows
$$\dim\Fi{K_{S_0}+\roundup{\pi^*(K_X)|_{S_0}}+L_0}(C)=1,$$
thus
$\dim\Fi{G}(C)=1$. Therefore
$\dim\Fi{M_{2k_0+2}|_{S_0}}(C)=1.$
and so $M_{2k_0+2}|_{S_0}\cdot C\ge 2$.
Noting that
$h^0(S_0, M_{2k_0+2}|_{S_0})\ge h^0(S_0,G)\ge 2,$
we get from Lemma 1.4 that
$K_{S_0}+\roundup{t\pi^*(K_X)|_{S_0}}+M_{2k_0+2}|_{S_0}$
is effective.
Now Lemma 1.3 implies the birationality of the rational map defined by
$|K_{S_0}+\roundup{t\pi^*(K_X)|_{S_0}}+M_{2k_0+2}|_{S_0}+L_0|.$
Because
\begin{align*}
&|K_{S_0}+\roundup{t\pi^*(K_X)|_{S_0}}+M_{2k_0+2}|_{S_0}+L_0|\\
\subset &|K_{S_0}+\roundup{(t+2k_0+2)\pi^*(K_X)}|_{S_0}+L_0|,
\end{align*}
$\Fi{K_{S_0}+\roundup{(t+2k_0+2)\pi^*(K_X)}|_{S_0}+S_0|_{S_0}}$
is birational.
We have proved the theorem.
\end{proof}

From now on, we suppose that $\dim\fei{k_0}(X)=1$. We can take the
same modification $\pi:X'\lrw X$ as in the proof of Proposition
2.1. Let $g:=\fei{k_0}\circ\pi$ be the morphism from $X'$  onto
$W\subset {\Bbb P}^{P_{k_0}-1},$ where $W$ is the Zariski closure
of the image of $X$ through $\fei{k_0}$. Let $g:X'\overset{f}\lrw
Q\lrw W$ be the Stein-factorization. Then $Q$ is a smooth
projective curve. Denote $b:=g(Q)$, the genus of $Q$.

\begin{setup}\label{Set:2.5}
If $b>0$, we have already known from \cite{Ch2} that $k_s\le 2k_0+4$.
\end{setup}

In the rest of this section, we mainly study the case when $Q$ is the
rational curve ${\Bbb P}^1$. We have a derived fibration $f:X'\lrw{\Bbb P}^1$. Let
$S$ be a general fiber of the fibration. Then $S$  is a smooth projective
surface of general type. Note that $S$ is also the generic irreducible
element of the movable part of the system $|k_0K_{X'}|$.
Let $\sigma: S\lrw S_0$ be the contraction onto the minimal model.

\begin{thm}\label{T:2.6} Let $X$ be a minimal projective 3-fold of general
type with only ${\Bbb Q}$-factorial terminal singularities.
If $\dim\fei{k_0}(X)=1$ and $b=0$, then $k_s\le 5k_0+6$.
\end{thm}
\begin{proof}
For all $i>0$, denote by $M_i$ the movable part of $|iK_{X'}|$. By
Koll\'ar's method (\cite{Kol} or see \cite{Ch2}, 2.2), we have
$M_{9k_0+4}|_S\ge 4\sigma^*(K_{S_0})$.

By the vanishing theorem, one has
\begin{align*}
&|K_{X'}+\roundup{(9k_0+4)\pi^*(K_X)}+S||_S=
|K_S+\roundup{(9k_0+4)\pi^*(K_X)}_S|\\
&\supset |K_S+M_{9k_0+4}|_S|\supset
|5\sigma^*(K_{S_0})|.
\end{align*}
By Lemma 1.5, we see that $M_{10k_0+5}|_S\ge
5\sigma^*(K_{S_0}).$ Repeatedly performing this process, one has
$$  M_{9k_0+4+m(k_0+1)}|_S\ge (m+4)\sigma^*(K_{S_0})$$
for all integer $m>0$.
This means that we can write
$$(5k_0+(m+4)(k_0+1))\pi^*(K_X)|_S\sim_{\Bbb Q} (m+4)\sigma^*(K_{S_0})+
E_{\Bbb Q}^{(m)},$$
where $E_{\Bbb Q}^{(m)}$ is an effective ${\Bbb Q}$-divisor only relating to
$m$. Thus
$$(\frac{5k_0}{m+4}+(k_0+1))\pi^*(K_X)|_S\simnum \sigma^*(K_{S_0})+
\frac{1}{m+4}E_{\Bbb Q}^{(m)}. \eqno(2.6.1)$$
We can write
$2\sigma^*(K_{S_0})\simlin C+Z,$
where $C$ is the movable part and $Z$ the fixed one. According to
Xiao(\cite{X}), $|C|$ is composed of a pencil of curves if and only if
$K_{S_0}^2=1$ and $p_g(S)=0$.
We prove this theorem step by step.
\smallskip

Step 1. $tK_{X'}$ is effective for all $t\ge 3k_0+4$.
\smallskip

Denote by $\alpha\ge 2k_0+3$ a positive integer. We have
$$|(\alpha+k_0+1)K_{X'}|\supset |K_{X'}+\roundup{\alpha\pi^*(K_X)}+S|.$$
By the vanishing theorem, one has
$$|K_{X'}+\roundup{\alpha\pi^*(K_X)}+S||_S=
|K_S+\roundup{\alpha\pi^*(K_X)}|_S|\supset |K_S+\roundup{\alpha\pi^*(K_X)|_S}|.$$
Now we consider a sub-system
$$|K_{S}+\roundup{\alpha\pi^*(K_X)|_S-Z-\frac{2}{m+4}E_{\Bbb Q}^{(m)}}|.$$
Because
$$\alpha\pi^*(K_X)|_S-Z-\frac{2}{m+4}E_{\Bbb Q}^{(m)}-C\simnum
t_0\pi^*(K_X)|_S$$
where $t_0:=\alpha-\frac{10k_0}{m+4}-2k_0-2>0$ for big $m$. Thus we have
$$|K_{S}+\roundup{\alpha\pi^*(K_X)|_S-Z-\frac{2}{m+4}E_{\Bbb Q}^{(m)}}||_S=
|K_C+D|,$$
where $D$ is divisor on $C$ with $\deg(D)>0$. Therefore
$P_{\alpha+k_0+1}(X)\ge h^0(C, K_C+D)\ge 2$ since $g(C)\ge 2$.
\smallskip

Step 2. The birationality.
\smallskip

Denote by $\beta\ge 4k_0+5$ a positive integer. Considering the system
$|K_{X'}+\roundup{\beta\pi^*(K_X)}+S|$, we have
$$ |K_{X'}+\roundup{\beta\pi^*(K_X)}+S||_S\supset
|K_S+\roundup{\beta\pi^*(K_X)|_S}|.$$
By Step 1, it's sufficient to verify the birationality of
$\Phi_{K_S+\roundup{\beta\pi^*(K_X)|_S}}$.
It's obvious that
$$K_S+\roundup{\beta\pi^*(K_X)|_S}\ge K_S+
\roundup{\beta\pi^*(K_X)|_S
-Z-\frac{4}{m+4}E_{\Bbb Q}^{(m)}}.$$
Denote by $A:=
\roundup{\beta\pi^*(K_X)|_S
-Z-\frac{4}{m+4}E_{\Bbb Q}^{(m)}}-2\sigma^*(K_{S_0})$. We have
$$A\simnum (\beta-4k_0-4-\frac{20k_0}{m+4})\pi^*(K_X)|_S.$$
So $A$ is also nef and big for big $m$. Now we have
$$|K_S+\roundup{\beta\pi^*(K_X)|_S-Z-\frac{4}{m+4}E_{\Bbb Q}^{(m)}}|=
|K_S+\roundup{A}+2\sigma^*(K_{S_0})+C|.$$
By Lemmas 1.3 and 1.4, one can easily see that the above system defines a birational
map onto its image. The theorem follows.
\end{proof}

\section{\rm The case with $k_0=1$}\label{S:3}

{}From now on, we only consider the case with $p_g(X)\ge 2$. We always suppose $X$ is a minimal
projective 3-fold of general type with ${\Bbb Q}$-factorial terminal singularities.
The first effective result was obtained  by Koll\'ar  who proved that $\fei{16}$
is birational. In \cite{Ch3}, it has been proved that $\fei{9}$ is birational. Here
we would like to prove the birationality of $\fei{8}$.

By virtue of Theorem 2.3, Theorem 2.4 and (2.5), we only have to consider the situation
with $\dim\phi_1(X)=1$
and $b=0$.  We have a derived fibration $f:X'\lrw C$.  A general fiber $S$ is a
projective smooth surface of general type. We note that $p_g(S)>0$ in this case.
In order to formulate our proof, we classify $S$ into the following types:

(I) $K_{S_0}^2=1$, $p_g(S)=2$;

(II) $K_{S_0}^2=2$, $p_g(S)=3$;

(III) $K_{S_0}^2=2$, $p_g(S)=2$;

(IV) $K_{S_0}^2=1$, $p_g(S)=1$;

(V) $K_{S_0}^2\ge 3$;

(VI)  $K_{S_0}^2=2$, $p_g(S)=1$.

\begin{claim}\label{Claim:3.1} If $S$ is of type (I), then $\Fi{8K_X}$ is
birational.
\end{claim}
\begin{proof}
Denote by $C$ a generic irreducible element of the movable part of $|K_S|$. Then it's well-known
that $C$ is a smooth curve of genus 2.  According to Claim in Proposition 5.3 of \cite{Ch3}, we have
$\xi:=\pi^*(K_X)\cdot C\ge \frac{3}{5}$. For a positive integer $t$, we have
$$|K_{X'}+\roundup{t\pi^*(K_X)}+S||_S=|K_S+\roundup{t\pi^*(K_X)}|_S|\supset
|K_S+\roundup{t\pi^*(K_X)|_S}|.$$
Since $p_g(X)>0$, taking $t=1$ and applying Lemma 1.5, one has $M_3|_S\ge C$. Taking $t=3$
and applying Lemma
1.5 once more, one has $M_5|_S\ge 2C$. This means
$$5\pi^*(K_X)|_S\sim_{\Bbb Q}2C+E_{\Bbb Q}^{(5)},$$
where $E_{\Bbb Q}^{(5)}$ is an effective ${\Bbb Q}$-divisor. Thus we have
$$\frac{5}{2}\pi^*(K_X)|_S\simnum C+\frac{1}{2}E_{\Bbb Q}^{(5)}.$$
Now $t\pi^*(K_X)|_S-C-\frac{1}{2}E_{\Bbb Q}^{(5)}\simnum (t-\frac{5}{2})\pi^*(K_X)|_S$. It is easy to see that, for $t\ge 6$,
$(t-\frac{5}{2})\pi^*(K_X)|_S\cdot C>2$. Lemma 1.3 implies that $\Fi{K_S+
\roundup{6\pi^*(K_X)|_S}}$ is birational and so is
$\Fi{8K_X}$.
\end{proof}

\begin{claim}\label{Claim:3.2} If $S$ is of type (II), then $\Fi{7K_X}$ is
birational.
\end{claim}
\begin{proof}
We still denote by $C$ a generic irreducible element of the movable part of $|K_S|$. It's well-known that
$|C|$ defines a generically finite map and $C$ is a smooth curve of genus 3. By a parallel argument as in
the proof of Claim 3.1, we have $M_{2m+1}|_S\ge mC$ for any positive integer $m$. This means
$$(2m+1)\pi^*(K_X)|_S\sim_{\Bbb Q}mC+E_{\Bbb Q}^{(m)},$$
where  $E_{\Bbb Q}^{(m)}$ is an effective ${\Bbb Q}$-divisor depending on $m$. Thus we have
$$\frac{2m+1}{m}\pi^*(K_X)|_S\simnum C+\frac{1}{m}E_{\Bbb Q}^{(m)}.$$
Therefore $\eta:=\pi^*(K_X)|_S\cdot C\ge \frac{m}{2m+1}C^2\ge \frac{2m}{2m+1}$ for all $m>0$.
So $\eta\ge 1$.
We want to verify the birationality of $|K_S+\roundup{t\pi^*(K_X)|_S}|$ for certain $t$. Because
$$t\pi^*(K_X)|_S-C-\frac{1}{m}E_{\Bbb Q}^{(m)}\simnum (t-\frac{2m+1}{m})\pi^*(K_X)|_S=:A,$$
Fix a big $m$, one can see that $A\cdot C>2$ for $t\ge 5$. Thus, by Lemma 1.3,
$\Fi{K_S+\roundup{5\pi^*(K_X)|_S}}$ is birational and so is $\Fi{7K_X}$.
\end{proof}

\begin{claim}\label{Claim:3.3}  If $S$ is of type (III), then $\Fi{7K_X}$ is
birational.
\end{claim}
\begin{proof}
Denote by $C$ a generic irreducible element of the movable part of $|K_S|$.
Recall that $\sigma:S\lrw S_0$ is the contraction onto the minimal model. $C_1:=\sigma_{*}(C)$ is the movable part of
$|K_{S_0}|$. It's easy to see that $C_1$ has two types:

(3.3.1) $|K_{S_0}|=|C_1|+Z$, where $C_1^2=0$ and $C_1$ is smooth curve of genus 2.

(3.3.2) $|K_{S_0}|=|C_1|$ , where $C_1$ is a smooth curve of genus 3.

In either cases, we always have $\sigma^*(K_{S_0})\cdot C=K_{S_0}\cdot C_1=2$. By
Theorem 3.1 of \cite{Ci},
$|mK_{S_0}|$ is basepoint free for $m\ge 2$. Now Koll\'ar's technique gives
$7\pi^*(K_X)|_S\ge 2\sigma^*(K_{S_0})$ and so $M_7|_S\ge 2\sigma^*(K_{S_0})$. By a parallel
argument as in the proof of Claim 3.1, we get
$$(2m+3)\pi^*(K_X)|_S\ge M_{2m+3}|_S\ge m\sigma^*(K_{S_0})$$
for all positive integer $m$. Thus
$$\pi^*(K_X)|_S\cdot C\ge \frac{m}{2m+3}\sigma^*(K_{S_0})\cdot C\ge \frac{2m}{2m+3}.$$
So $\pi^*(K_{S_0})\cdot C\ge 1$. Now by the same argument as in the proof of Claim 4.2, one
can easily obtain the birationality of $\Fi{7K_X}$.
\end{proof}

\begin{claim}\label{Claim:3.4} If $S$ is of type (IV), then $\Fi{8K_X}$ is
birational.
\end{claim}
\begin{proof}
We consider the natural map
$$H^0(X',M_2)\overset{\alpha_2}{\lrw}\Lambda_2\subset H^0(S,M_2|_S)\subset
H^0(S,2K_S),$$
 where $\Lambda_2$ is the image of $\alpha_2$. Since $P_2(S)=3$,
$0<\dim_{\Bbb C}\Lambda_2\le 3$.

(3.4.1) $\dim_{\Bbb C}\Lambda_2=3$. In this case, $\Lambda_2$ defines the bicanonical map of
$S$. By \cite{Ci}, $|2K_{S_0}|$ is base point free. Thus we see that $M_2|_S\ge 2\sigma^*(K_{S_0})$.
We can write $2\pi^*(K_X)|_S\sim 2\sigma^*(K_{S_0})+E_{\Bbb Q}$, where $E_{\Bbb Q}$ is an effective
${\Bbb Q}$-divisor. Denote by $C$ a general member of $|2\sigma^*(K_{S_0})|$.
Now we have $4\pi^*(K_X)|_S-C-E_{\Bbb Q}\simnum 2\pi^*(K_X)|_S$ and
$2\pi^*(K_X)|_S\cdot C\ge 4$. By Lemma 1.3, $\Fi{K_S+\roundup{4\pi^*(K_X)|_S}}$ is birational and so is
$\Fi{6K_X}$.

(3.4.2) $\dim_{\Bbb C}\Lambda_2=2$. Because $\Lambda_2$ defines a morphism,
the movable part of $\Lambda_2$ forms a complete linear pencil. The pencil is rational because
$q(S)=0$. Denote by $C_1$ a generic irreducible element of the movable part of
$\Lambda_2$. By Lemma 3.7 below, $\sigma^*(K_{S_0})\cdot C_1\ge 2$.
Because $M_2|_S\ge C_1$, we can write $2\pi^*(K_X)|_S\simnum C_1+E_{\Bbb Q}'$, where
$E_{\Bbb Q}'$ is an effective ${\Bbb Q}$-divisor.
Now $5\pi^*(K_X)|_S-C_1-E_{\Bbb Q}'\simnum 3\pi^*(K_X)|_S$ and, by
(3.6.1), $\pi^*(K_X)|_S\cdot C_1\ge \frac{1}{2}\sigma^*(K_{S_0})\cdot C_1\ge 1$.
According to Lemma 1.3,  $\Fi{K_S+\roundup{5\pi^*(K_X)|_S}}$ is birational and so is
$\Fi{7K_X}$.

(3.4.3) $\dim_{\Bbb C}\Lambda_2=1$. In this case, $|2K_X|$ is composed of
a pencil of surfaces. One can see that $q(X)=h^2(\Co{X})=0$, whence
$\chi(\Co{X})\le -1$. Applying Reid's R-R formula(\cite{R2}), one has
$P_2(X)\ge 4$.
We can write $2\pi^*(K_X)\simnum aS+E'$,where $a\ge 3$ and $E'$ is an effective divisor.
We hope to prove the birationality of $\Fi{8K_X}$. Because $p_g(X)>0$,
it's sufficient to prove the birationality of $\Fi{8K_X}|_S$. By virtue of
Koll\'ar's method, one can see that $|7K_{X'}||_S\supset |2\sigma^*(K_{S_0})|$.
Therefore we are reduced to prove the birationality of $(\Fi{8K_{X'}}|_S)|_C$, where
$C$ is a general member of $|2\sigma^*(K_{S_0})|$.

By the vanishing theorem, one has
$|K_{X'}+\roundup{7\pi^*(K_X)-\frac{1}{a}E'}||_S\supset
|K_S+\roundup{L}|$ where $L\simnum (7-\frac{2}{a})\pi^*(K_X)|_S$.
(3.5.1) below is still true when $S$ is of type (IV),  i.e.
$$\frac{2(2m+3)}{m}\pi^*(K_X)\simnum C+E_{\Bbb Q}^{(m)},$$
where $E_{\Bbb Q}^{(m)}$ is an effective ${\Bbb Q}$-divisor.
Thus $|K_S+\roundup{L-E_{\Bbb Q}^{(m)}}||_C=|K_C+D|$, where
$D:=\roundup{L-E_{\Bbb Q}^{(m)}}|_C$. Because
$L-C-E_{\Bbb Q}^{(m)}\simnum (3-\frac{2}{3}-\frac{6}{m})\pi^*(K_X)|_S$ for all $m>0$,
we see that $\deg(D)>2$ whenever $m$ is large. Thus $\Fi{K_S+\roundup{L}}|_C$ is birational and so is
$\Fi{8K_X}$.
\end{proof}

\begin{claim}\label{Claim:3.5} If $S$ is of type (V), then $\Fi{7K_X}$ is
birational.
\end{claim}
\begin{proof}
By \cite{Ci}, $|mK_{S_0}|$ is basepoint free for all $m\ge 2$.
Denote by $C$ the movable part of $|2K_S|$. Then $C\simlin
2\sigma^*(K_{S_0})$. It's sufficient to prove the birationality of
$\Fi{7K_X}|_S$. By Koll\'ar's method, $|7K_{X'}||_S\supset
|2\sigma^*(K_{S_0})|=|C|$. We only need to verify the
birationality of $(\Fi{7K_{X'}}|_S)|_C $. The vanishing theorem
gives
$|K_{X'}+\roundup{7\pi^*(K_X)}+S||_S\supset|K_S+\roundup{7\pi^*(K_X)|_S}|.$
Applying Lemma 1.5, we get $M_9|_S\ge 3\sigma^*(K_{S_0})$.
Repeatedly proceeding the above process while replacing "7" by "9,
11, $\cdots$",  one can obtain $M_{2m+3}|_S\ge m\sigma^*(K_{S_0})$
and so
$$(2m+3)\pi^*(K_X)|_S\simnum m\sigma^*(K_{S_0})+E_{\Bbb Q}^{(m)},$$
where $m$ is a positive integer. Thus we have
$$\frac{2(2m+3)}{m}\pi^*(K_X)|_S\simnum 2\sigma^*(K_{S_0})+
\frac{2}{m}E_{\Bbb Q}^{(m)}.\eqno(3.5.1)$$
It's easy to see that $\pi^*(K_X)|_S\cdot C\ge 3$.

Now taking a very big $m$, one has
$$|K_S+\roundup{5\pi^*(K_X)|_S-\frac{2}{m}E_{\Bbb Q}^{(m)}}||_C=|K_C+D|,$$
where $D:=\roundup{5\pi^*(K_X)|_S-\frac{2}{m}E_{\Bbb Q}^{(m)}}|_C$  and
$\deg(D)\ge (1-\frac{6}{m})\pi^*(K_X)|_S\cdot C>2.$
Therefore we have proved that $(\Fi{7K_{X'}}|_S)|_C$ is birational. So $\Fi{7K_X}$ is birational.
\end{proof}

\begin{claim}\label{Claim:3.6} If $S$ is of type (VI), then $\Fi{8K_X}$ is
birational.
\end{claim}
\begin{proof}
We keep the same notations as in the proof of Claim 3.5. The proof is almost the same except
that we have here $\pi^*(K_X)|_S\cdot C\ge 2$. Thus we can prove that
$|K_S+\roundup{6\pi^*(K_X)|_S-\frac{2}{m}E_{\Bbb Q}^{(m)}}||_C$ is birational by the same
argument.  This, in turn,  proves the birationality of $\Fi{8K_X}$. We conclude the claim.
\end{proof}

Claims 3.1 through 3.6 imply Theorem 0.2.

\begin{lemma}\label{L:3.7}
Let $S$ be a smooth projective surface of general type. Let $\sigma: S\lrw S_0$ be the
contraction onto the minimal model. Suppose we have an effective irreducible curve $C$
on $S$ such that $C\le \sigma^*(2K_{S_0})$ and $h^0(S, C)=2$. If $K_{S_0}^2=p_g(S)=1$,
then $C\cdot \sigma^*(K_{S_0})\ge 2$.
\end{lemma}
\begin{proof}
We can suppose $|C|$ is a free pencil. Otherwise, we can blow-up $S$ at base points of
$|C|$. Denote $C_1:=\sigma(C)$. Then $h^0(S_0, C_1)\ge 2$. Suppose
$C\cdot \sigma^*(K_{S_0})=1$. Then $C_1\cdot K_{S_0}=1$. Because $p_a(C_1)\ge 2$,
we can see that $C_1^2>0$. {}From $K_{S_0}(K_{S_0}-C_1)=0$, we get
 $(K_{S_0}-C_1)^2\le 0$, i.e. $C_1^2\le 1$. Thus $C_1^2=1$ and $K_{S_0}\simnum C_1$.
This means $K_{S_0}\simlin C_1$ by virtue of \cite{Ca}, which is impossible because $p_g(S)=1$.
So $C\cdot\sigma^*(K_{S_0})\ge 2$.
\end{proof}

\begin{remark}\label{R:3.8} Slightly modifying our method, one can even prove the following
statements:

Let $X$ be a minimal projective 3-fold of general type with ${\Bbb Q}$-factorial terminal singularities.
Suppose $p_g(X)\ge 2$. Then

(1) $\fei{7}$ is birational whenever $p_g(X)\ne 2$.

(2) $\fei{6}$ is birational whenever $p_g(X)\ne 2,\ 3$.

(3) $\fei{5}$ is birational whenever $p_g(X)\ne 2,\ 3,\  4$.

We don't give the explicit proof since we think that the calculation is more complicated.  However
it's not difficult for a reader to verify these statements once he understands our method.
\end{remark}

\centerline{\bf Acknowledgement}
 This paper was begun while I was doing a post-doc
 research at the Mathematical Institute of the University of Goettingen, Germany between
 1999 and 2000. It was finally
 revised when I visited the IMS, CUHK in 2003.  I appreciate very
 much for financial supports from both the institutes.
Thanks are also due to both Professor Fabrizio Catanese for
effective discussions and Professor Kang Zuo for many helps and
hospitality.

\end{document}